\documentclass[11pt]{article}
\usepackage[utf8]{inputenc}
\usepackage[T1]{fontenc}
\usepackage{amsmath}
\usepackage{amsfonts}
\usepackage{amssymb}
\usepackage[version=4]{mhchem}
\usepackage{stmaryrd}
\usepackage{graphicx}
\usepackage[export]{adjustbox}
\graphicspath{ {./images/} }
\usepackage{hyperref}
\hypersetup{colorlinks=true, linkcolor=blue, filecolor=magenta, urlcolor=cyan,}
\usepackage{bbold}
\usepackage{multicol}

\newtheorem{prop}{Proposition}
\newtheorem{cor}{Corollary}

\begin{document}

\begin{centering}
\Large
{YET ANOTHER SPECIES\\
OF FORBIDDEN-DISTANCES \\ 
CHROMATIC NUMBER}\footnote{This article was published in \emph{Geombinatorics}, vol.~10 no.~3 (2001), pp.~89--95.}

\large

\bigskip

\begin{multicols}{2}
{Aaron Abrams\\
Department of Mathematics\\ University of Georgia \\Athens, Georgia
\\{\tt abrams.aaron@gmail.com}}

Peter Johnson\\
Department of Discrete and Statistical Sciences \\
Auburn University, AL 36849 
{\tt johnspd@mail.auburn.edu} 

\end{multicols}

\end{centering}

\normalsize

\bigskip
\bigskip

Suppose that $(X, \rho)$ is a metric space. We will follow the custom of sometimes suppressing mention of the metric $\rho$, and letting $X$ refer to the pair $(X, \rho)$. If $s>0$, the chromatic number of $X$ with respect to $s$, denoted $\chi(X, s)$, is the smallest number of colors required to color $X$ so that no two points of $X$ a distance $s$ apart are the same color. We then define the chromatic number of $X$ to be $\chi(X)=\sup \chi(X, s)$. [Every "sup" in this paper is a "max", if finite.] When $X=\mathbb{R}^{n}$ and $\rho$ is a metric arising from a norm, all the values $\chi(X, s), s>0$, are the sane, so our definition agrees with common usage in these cases.

The \emph{lower chromatic number} of $X$ is the smallest positive integer $m$, if any, such that there exist distances $s_{1}, \ldots, s_{m}>0$, actually realized between points of $X$, and a partition (coloring) of $X$ into sets $X_{1}, \ldots, X_{m}$ such that no two points of $X_{i}$ are a distance $s_{i}$ apart, $i=1, \ldots, m$. (We could omit the provision ``actually realized between points of $X$,'' but then the lower chromatic number of any metric space in which not all distances are realized would be 1.) If no such $m$ exists, we declare the lower chromatic number of $X$ to be $\infty$. Letting $\underline \chi$ denote the lower chromatic number, clearly we have $\underline \chi(X) \leq \chi(X)$, if $|X| \geq$ 2. In the case $X=\mathbb{R}^{2}, \underline{\chi}$ is the same as the \emph{polychromatic number}, in [6].

The upper chromatic number of $X$, denoted $\hat{\chi}(X)$, is the smallest positive integer $m$, if any, such that for \emph{every} sequence $s_{1}, \ldots, s_{m}>0$, there is a coloring of $X$ with $m$ colors such that the distance $s_{i}$ is forbidden for the $i^ {\text{th}}$ color, $i=1, \ldots, m$; to put it another way, for each sequence $s_{1}, \ldots, s_{m}>0$ there is a partition of $X$ into sets $X_{1}, \ldots, X_{m}$ (allowed to be empty) such that no two points of $X_{i}$ are a distance $s_{i}$ apart, $i=1, \ldots, m$. If no such $m$ exists, we set $\hat{\chi}(X)=\infty$. Clearly $\chi(X) \leq \hat{\chi}(X)$ : supposing $m=\hat{\chi}(X)$ is finite, given $s>0$, considering the sequence $s=s_{1}=\cdots=s_{m}$, we see that $X$ is $m$-colorable so that the distance $s$ is forbidden for all colors; thus $\chi(X, s) \leq$ $\hat{\chi}(X)$ for all $s>0$.

More generally, for $k=1,2, \ldots$, the $k^{\text{th}}$ upper chromatic number of $X$, $\hat{\chi}^{(k)}(X)$, is defined to be the smallest positive integer $m$ (as before, $\hat{\chi}^{(k)}(X)=\infty$ if there is no such $m$) such that for every $k \times m$ matrix of positive numbers, there is a coloring of $X$ with $m$ colors such that for each $i$, the distances in the $i^{\text{th}}$ column of the matrix are forbidden for the $i ^{\text{th }}$ color, $i=1, \ldots, m$. Clearly $\hat{\chi}^{(k)}(X) \leq \hat{\chi}^{(k+1)}(X), k=1,2, \ldots$, and $\hat{\chi}^{(1)}=\hat{\chi}$.

Only four papers, [1] - [4], have been written so far about the parameters $\hat{\chi}^{(k)}$; [2] will serve as a survey and starting point, for those interested. The main results on the $\hat{\chi}^{(k)}$ are that $\hat{\chi}^{(k)}(\mathbb{R})$ is finite for all $k$ (but little is known about the actual values, beyond $\hat{\chi}(\mathbb{R})=3,\ 4 \leq \hat{\chi}^{(2)}(\mathbb{R}) \leq 22$, and $k+1 \leq \chi^{(k)}(\mathbb{R}) \leq$ $\lceil 4 e k\rceil)$, $\hat{\chi}^{(k)}\left(\mathbb{R}^{n}\right)$ is finite for all $n$ and $k$ if the metric on $\mathbb{R}^{n}$ is the $\ell^{\infty}$ metric, arising from the max norm; but it is not even known, for instance, whether or not $\hat{\chi}\left(\mathbb{Z}^{2}\right)$ is finite, if $\mathbb{Z}^{2}$ bears the usual Euclidean metric.

Our purpose here is to promote the study of another sequence of ``forbidden-distances'' chromatic numbers related to the upper chromatic numbers, but evidently much more tractable than they are. Already introduced in [2], the (quite natural) definition of these parameters is due to Babai, in a communication with the first author. For a positive integer $k$ and a metric space $(X, \rho), B_{k}(X)$ is the smallest positive integer among those $m$, if any, such that for every set of $k$ or fewer positive numbers, there is a coloring of $X$ with $m$ colors with every distance in the set forbidden for every color. If there is no such $m$, we declare $B_{k}(X)=\infty$.

Obviously $B_{k}(X) \leq B_{k+1}(X)$ and if $X \subseteq Y$, with the metric on $X$ being the restriction to $X \times X$ of the metric on $Y$, then $B_{k}(X) \leq B_{k}(Y), k=1,2, \ldots$. Further, $B_{1}(X)=$ $\chi(X)$. Indeed, $B_{k}(X)$ is the supremum of the chromatic numbers $\chi\left(X ; s_{1}, \ldots, s_{k}\right)$ of the graphs obtained by taking the points of $X$ as vertices, with any two vertices adjacent if and only if the distance between them is one of the numbers $s_{1}, \ldots, s_{k}$. The supremum is taken over all choices of $s_{1}, \ldots, s_{k}>0$. It is worth noting that neither $\underline{\chi}$ nor the $\hat{\chi}^{(k)}$ can be viewed (at least not in any obvious way) as being related to ordinary chromatic numbers of graphs, as are the $B_{k}$.

\begin{prop}
For $k=1,2, \ldots,\ \chi(X) \leq B_{k}(X) \leq \hat{\chi}^{(k)}(X)$.
\end{prop}

{\bf Proof:} By previous remarks, $\chi(X)=B_{1}(X) \leq B_{k}(X)$. If $\hat{\chi}^{(k)}(X)=m<\infty$ and $s_{1}, \ldots, s_{k}>0$, consider the $k \times m$ matrix with $s_{1}, \ldots, s_{k}$ in each column. There is an $m$-coloring of $X$ with the distances in the $j^{\text {th}}$ column of the matrix forbidden for the $j^{\text{th}}$ color; so the distances $s_{1}, \ldots, s_{k}$ are forbidden for every color. This proves that $B_{k}(X) \leq m=\hat{\chi}^{(k)}(X)$. $\square$

\begin{prop} 
For positive integers $k, p$, $B_{k+p}(X) \leq B_{k}(X) B_{p}(X)$.
\end{prop}

{\bf Proof:} Suppose $B_{k}(X)=m$ and $B_{p}(X)=r$ are both finite. Suppose $s_{1}, \ldots, s_{k}, s_{k+1}, \ldots, s_{k+p}>0$. Color $X$ with $m$ colors so that the distances $s_{1}, \ldots, s_{k}$ are forbidden for each color. Let $X_{1}, \ldots X_{m}$ be the corresponding color classes $\left(X_{j}\right.$ is the set of points bearing the $j^{\text{th}}$ color): since $B_{p}\left(X_{j}\right) \leq$ $B_{p}(X)$, we can color each $X_{j}$ with $r$ colors so that the distances $s_{k+1}, \ldots, s_{k+p}$ are forbidden for each color. Let $X_{1 j}, \ldots, X_{r j}$ be the color classes within $X_{j}$; we now see that $X$ is partitioned into sets $X_{i j}, i=1, \ldots, r,\ j=1, \ldots, m$, with none of distances $s_{1}, \ldots, s_{k+p}$ occurring between points of any $X_{i j}$. $\square$

\begin{cor}
$B_{k}(X) \leq \chi(X)^{k}, k=1,2, \ldots$ .
\end{cor}

\begin{cor} 
$B_{k}\left(\mathbb{R}^{n}\right)$, with $\mathbb{R}^{n}$ bearing the Euclidean metric (or any metric such that $\left.\chi\left(\mathbb{R}^{n}\right)<\infty\right)$, is finite.
\end{cor}

In all that follows, the metric on $\mathbb{Z}^{n}, \mathbb{Q}^{n}$, or $\mathbb{R}^{n}$ will be the usual Euclidean metric, in each case.

\begin{prop}
$k+1 \leq B_{k}(\mathbb{Z}) \leq B_{k}(\mathbb{R}) \leq 2 k$ for $k=1,2, \ldots$ .
\end{prop}

{\bf Proof:} Trying to color $\{0,1, \ldots, k\}$ so as to forbid distances $1, \ldots, k$ shows that $k+1 \leq B_{k}(\mathbb{Z})$. If $s_{1}, \ldots, s_{k}>0$ are distances, and $G$ is the graph with vertices $\mathbb{R}$, with any two points $a,b$ adjacent if and only if $|a-b| \in\left\{s_{1}, \ldots, s_{k}\right\}$, then $\Delta(G) \leq 2 k$, so the same holds for any finite subgraph of $G$, and it is clear that $G$ contains no $K_{2 k+1}$ nor, in the case $k=1$, any cycle at all, much less an odd cycle. Brooks' Theorem now implies that the chromatic number of any finite subgraph of $G$ is at most $2 k$, so the same is true of $G$, by a famous theorem of de Bruijn and Erd\"os. Since $s_{1}, \ldots, s_{k}$ were arbitrary, it follows that $B_{k}(\mathbb{R}) \leq 2 k$. $\square$

\bigskip

The lower bound in Proposition 3 is a special case of a more general lower bound, introduced in [2] to apply to the upper chromatic numbers but which more naturally, and perhaps more sharply, applies to the $B_{k}$. A subset $S$ of a metric space $X$ is a \emph{$k$-distance set} if and only if there are no more than $k$ non-zero distances among the points of $S$, i.e., $|\{\rho(x, y)\ ;\ x, y \in S, x \neq y\} | \leq k$. So, for instance, $\{0, \ldots, k\}$ is a $k$-distance set in $\mathbb{Z}$ (with the usual metric). (Note that there could be no larger $k$-distance set in $\mathbb{R}$.)

\begin{prop}
$B_{k}(X) \geq \sup \{\ |S|\ ;\ S \subseteq X$ is a k-distance set $\}$. 
\end{prop}
The proof is clear, and is omitted.

\begin{cor}
$B_{k}\left(\mathbb{R}^{2}\right) \geq 2 k+1 ;\ B_{k}\left(\mathbb{R}^{k}\right) \geq 2^{k} ;\ B_{k}\left(\mathbb{R}^{n}\right) \geq {\binom{n+1} k} ;
\ B_{3}\left(\mathbb{R}^{3}\right) \geq$ 12.
\end{cor}

For proofs, see [2].

Finally, we give some values or estimates of $B_{k}(X)$ for $k \in$ $\{1,2\}$ and $X \in\left\{\mathbb{Z}, \mathbb{R}, \mathbb{Z}^{2}, \mathbb{Q}^{n}, \mathbb{Z}^{3}\right\}$.

\begin{prop}
(a) $B_{2}(\mathbb{Z})=B_{2}(\mathbb{R})=3$;

(b) $B_{1}\left(\mathbb{Z}^{2}\right)=B_{1}\left(\mathbb{Q}^{2}\right)=2<B_{1}\left(\mathbb{Z}^{3}\right)$;

(c) $B_{2}\left(\mathbb{Z}^{2}\right)=B_{2}\left(\mathbb{Q}^{2}\right)=4$.
\end{prop}

{\bf Proof:}
\begin{itemize}
\item[(a)] $3 \leq B_{2}(\mathbb{Z})$ by Proposition 3. Now suppose that $0<s_{1} \leq s_{2}$. We want to color $\mathbb{R}$ with 3 colors so that both distances $s_{1}, s_{2}$ are forbidden, for all colors. If $s_{1}=s_{2}$ we can do the coloring with two colors, so suppose that $s_{1}<s_{2}$. We have $s_{2}=m s_{1}+a$ for some positive integer $m$ and some $a \in\left[0, s_{1}\right)$.

Partition $\mathbb{R}$ into half-open intervals of length $s_{2}, J_{n}=\left[n s_{2},(n+1) s_{2}\right)$, $n \in \mathbb{Z}$. Partition each $J_{n}$ into half-open intervals, $m$ of length $s_{1}$ and a last one of length $a$, if $a>0$. Thus, these intervals are $$\left[n s_{2}+(q-1) s_{1}, n s_{2}+q s_{1}\right),$$ $q=1, \ldots, m$, and if $a>0$, the last is $\left[(n+1) s_{2}-a,(n+1) s_{2}\right)$. 

The three colors we will use will be red, blue, and green. Each of the small half-open intervals will be colored with one of these colors. Only two colors will appear in $J_{n}$. Starting on the left in $J_{n}$, and scanning to the right, color the small subintervals alternately red and blue (starting with red) if $n \equiv 0 \bmod 3$, green and red (starting with green) if $n \equiv 1 \bmod 3$, and blue and green (starting with blue) if $n \equiv 2 \bmod 3$. It is straightforward to check that no two points the same color are a distance $s_{1}$ or $s_{2}$ apart with $\mathbb{R}$ so colored. 

\item[(b)] Clearly $B_{1}\left(\mathbb{Z}^{2}\right) \geq 2$.

If $r \in \mathbb{Q},\ r>0$, and $\mathbb{Q}^{2}$ is two-colorable so as to forbid the distance $s>0$, then $\mathbb{Q}^{2}$ is two-colorable so as to forbid the distance $r s$. Therefore, to show that $B_{1}\left(\mathbb{Q}^{2}\right)=2$, it suffices to show that $\chi\left(\mathbb{Q}^{2}, s\right) \leq 2$ for two kinds of distance $s$:  $\sqrt{p / q}$ and $\sqrt{2 p / q},\ $ $p, q$ odd positive integers. The former is taken care of in [5]; in fact, it is possible to two-color $\mathbb{Q}^{2}$ (even $\mathbb{Q}^{3}$) so that all of those distances are forbidden at once.

Suppose that $p$ and $q$ are odd positive integers, and suppose that $\sqrt{2 p / q}$ occurs as a distance in $\mathbb{Q}^{2}$. Suppose that $a, b, c \in \mathbb{Z}, c>0$, with $\operatorname{gcd}(a, b, c)=1$, and $q\left(a^{2}+b^{2}\right)=$ $2 p c^{2}$. Since $\operatorname{gcd}(a, b , c)=1$, if $c$ is even then $a$ and $b$ cannot both be even, so $q\left(a^{2}+b^{2}\right) \not \equiv 0 \bmod 4$, while $2 p c^{2} \equiv 0 \bmod 4$; so $c$ must be odd. Therefore $q\left(a^{2}+b^{2}\right) \equiv 2 \bmod 4$; it must be that both $a$ and $b$ are odd.

Now suppose that $\left(\frac{a_{i}}{c_{i}}, \frac{b_{i}}{c_{i}}\right) \in \mathbb{Q}^{2}, \operatorname{gcd}\left(a_{i}, b_{i}, c_{i}\right)=1,\left(\frac{a_{i}}{c_{i}}\right)^{2}+$ $\left(\frac{b_{i}}{c_{i}}\right)^{2}=\frac{2 p}{q}, m_{i} \in \mathbb{Z}, i=1, \ldots, n$, and $\sum\limits_{i=1}^{n} m_{i}\left(\frac{a_{i}}{c_{i}}, \frac{b_{i}}{c_{i}}\right)=$ $(0,0)$. Then $\sum\limits_{i=1}^{n} m_{i} \frac{a_{i}}{c_{i}}=0$. Noting that all $a_{i}, c_{i}$ are odd, and multiplying the last equation through by the product of the $c_{i}$, we have

$$
0=\sum_{i=1}^{n} m_{i} \text { (odd integer) } \equiv \sum_{i=1}^{n} m_{i} \equiv \sum_{i=1}^{n}\left|m_{j}\right| \ \bmod 2
$$

All of this goes to show that for any closed walk in $\mathbb{Q}^{2}$ in steps of length $\sqrt{2 p / q}$, there must be an even number of steps. By a well known characterization of bipartite graphs, it follows that $\chi\left(\mathbb{Q}^{2}, \sqrt{2 p / q}\right)=2$. This completes the proof that $B_{1}\left(\mathbb{Z}^{2}\right)=B_{1}\left(\mathbb{Q}^{2}\right)=2$. To see that $B_{1}\left(\mathbb{Z}^{3}\right) \geq 3$, observe that $(0,0,0),(1,1,0) ,(1,0,1)$ are the vertices of an equilateral triangle in $\mathbb{Z}^{3}$. 

\item[(c)] $B_{2}\left(\mathbb{Z}^{2}\right) \geq 4$ because the vertices of a square constitute a 2-distance set. On the other hand, by Proposition 2, $B_{2}\left(\mathbb{Q}^{2}\right) \leq\left(B_{1}\left(\mathbb{Q}^{2}\right)\right)^{2}=4$.\\ \phantom{boo} \hfill $\square$

\end{itemize}

\bigskip

What's next? The discovery of $B_{1}\left(\mathbb{Q}^{3}\right), B_{1}\left(\mathbb{Q}^{4}\right)$, and $B_{3}(\mathbb{Z})$ may be within reach. Also, it seems plausible that $B_{k}(\mathbb{Z})=$ $B_{k}(\mathbb{R})$ for all $k$, in view of the definition of these numbers as the maximum chromatic numbers of certain kinds of graphs. It seems intuitively clear that those graphs on $\mathbb{R}$ that you cannot get on $\mathbb{Z}$ will be easier to color than those on $\mathbb{Z}$; but we have no clear proof, yet.

\vglue .5in

\Large
\noindent
{\bf References}
\normalsize
\bigskip

\begin{enumerate}

\item[[1]\!\!\!] A.~Abrams, The $k^{\text th}$ upper chromatic number of the line, \emph{Discrete Mathematics} 169 (1997), 157--162.

\item[[2]\!\!\!]  A.~Abrams, Upper chromatic numbers: an update, \emph{Geombinatorics} 10 (2000), 4--11.

\item[[3]\!\!\!]  A.~Archer, On the upper chromatic numbers of the reals, \emph{Discrete Mathematics} 214 (2000), 65--75.

\item[[4]\!\!\!]  D.~Greenwell and P.~Johnson, Forbidding prescribed distances for designated colors, \emph{Geombinatorics} 2 (1992), 13--16.

\item[[5]\!\!\!]  P.~Johnson, Two-colorings of a dense subgroup of $\mathbb{Q}^{n}$ that forbid many distances, \emph{Discrete Mathematics} 79 (1989/90), 191--195.

\item[[6]\!\!\!]  A.~Soifer, Relatives of the chromatic number of the plane I. \emph{Geombinatorics} 1 (1991), 13--17.

\end{enumerate}

\end{document}